\documentclass[12pt,a4paper]{article}

\usepackage{amsmath,amsfonts, amssymb}
\begin{document}
\bibliographystyle{abbrv}

\title{Preservation of log-concavity on summation}
\author{Oliver Johnson\thanks{Statistical Laboratory, 
Centre for Mathematical Sciences, University of Cambridge, Wilberforce Rd, 
Cambridge, CB3 0WB, UK.} 
\thanks{Christ's College, Cambridge. Email: {\tt otj1000@cam.ac.uk}\;
Fax: +44 1223 337956\; Phone: +44 1223 337946} 
\and Christina Goldschmidt$^*$\thanks{Pembroke College, Cambridge. Email: \tt{C.Goldschmidt@statslab.cam.ac.uk}}}
\date{\today}
\maketitle

\newtheorem{theorem}{Theorem}[section]
\newtheorem{lemma}[theorem]{Lemma}
\newtheorem{proposition}[theorem]{Proposition}
\newtheorem{corollary}[theorem]{Corollary}
\newtheorem{conjecture}[theorem]{Conjecture}
\newtheorem{definition}[theorem]{Definition}
\newtheorem{example}[theorem]{Example}
\newtheorem{condition}{Condition}
\newtheorem{main}{Theorem}
\newtheorem{remark}[theorem]{Remark}
\hfuzz20pt

\def \outlineby #1#2#3{\vbox{\hrule\hbox{\vrule\kern #1% 
\vbox{\kern #2 #3\kern #2}\kern #1\vrule}\hrule}}%
\def \endbox {\outlineby{4pt}{4pt}{}}%
\newenvironment{proof}
{\noindent{\bf Proof\ }}{{\hfill \endbox
}\par\vskip2\parsep}
\newenvironment{pfof}[2]{\removelastskip\vspace{6pt}\noindent
 {\it Proof  #1.}~\rm#2}{\par\vspace{6pt}}

\newcommand{\var}{{\rm{Var\;}}}
\newcommand{\cov}{{\rm{Cov\;}}}
\newcommand{\tends}{\rightarrow \infty}
\newcommand{\C}{{\cal C}}
\newcommand{\ep}{{\mathbb {E}}}
\newcommand{\pr}{{\mathbb {P}}}
\newcommand{\re}{{\mathbb {R}}}
\newcommand{\I}{\mathbb {I}}
\newcommand{\Z}{{\mathbb {Z}}}
\newcommand{\Ch}[2]{\ensuremath{\begin{pmatrix} #1 \\ #2 \end{pmatrix}}} 
\newcommand{\bin}[2]{\binom{#1}{#2}}
\newcommand{\ov}[1]{\overline{#1}}
\newcommand{\geo}[1]{{\rm{Geom}}\left(#1\right)}
\newcommand{\Prob}[1]{\ensuremath{\mathbb{P} \left(#1 \right)}}
\newcommand{\blah}[1]{}
\newcommand{\ch}[1]{{\bf #1}}

\begin{abstract}
\noindent
We extend Hoggar's theorem that the sum of two independent
discrete-valued log-concave random variables is itself log-concave. We
introduce conditions under which the result still holds for dependent
variables. We argue that these conditions are natural by giving some
applications.  Firstly, we use our main theorem to give simple proofs
of the log-concavity of the Stirling numbers of the second kind and of
the Eulerian numbers.
Secondly, we prove results concerning the log-concavity
of the sum of independent (not necessarily log-concave) random
variables.
\end{abstract}


\section{Introduction and main theorem}
The property of log-concavity (LC) for non-negative sequences 
is defined as follows.
\begin{definition} \label{def:logconc}
A sequence $(u(i), i \geq 0)$ is log-concave if, for all $i \geq 1$, 
\begin{equation} \label{eq:lcpropseq}
u(i)^2 \geq u(i-1) u(i+1).
\end{equation}
\end{definition}
Equation \ref{eq:lcpropseq} is sometimes referred to the quadratic
Newton inequality (see Niculescu \cite{niculescu}), and plays a role
in considering when polynomials with real coefficients have real roots
(see Section \ref{sec:appl}).
Applications of log-concavity arise in combinatorics, algebra,
geometry and computer science, as reviewed by Stanley~\cite{stanley}
and Brenti~\cite{brenti3,brenti2}.  In probability and statistics, it
is related to the notion of negative association of random variables
(see Joag-Dev and Proschan~\cite{joag-dev}).
\begin{definition} \label{def:logconcdisn}
A random variable $V$ taking values in $\Z_+$ is log-concave if its
probability mass function $p_V(i) = \pr(V=i)$ forms a log-concave sequence.
That is, $V$ is log-concave if for all $i \geq 1$, 
\begin{equation} \label{eq:lcpropdisn}
p_V(i)^2 \geq p_V(i-1) p_V(i+1).
\end{equation}
\end{definition}
\begin{example} \mbox{ } \label{ex:first}
\begin{enumerate}
\item 
We write $\geo{p}$ for the geometric distribution with mass function
$p_{X}(i) = (1-p) p^i$, for $i \in \Z_+$. For each $p$, these random
variables represent the `edge case' of the log-concavity property, in
that Equation (\ref{eq:lcpropdisn}) becomes an identity for all $i
\geq 1$.
\item The ${\rm Poisson}(\lambda)$ distribution is log-concave for any
$\lambda \geq 0$.
\item \label{ex:bernsum} Any Bernoulli random variable (that is, only
taking values in $\{ 0, 1 \}$) is log-concave. Further, any binomial
distribution is log-concave.  In fact, any random variable $S =
\sum_{i=1}^n X_i$, where $X_i$ are independent (not necessarily
identical) Bernoulli variables, is log-concave.
\end{enumerate}
\end{example}
This last fact is a corollary of the following theorem, first proved by Hoggar
\cite{hoggar}.
\begin{theorem} \label{thm:hoggar} If $V$ and $W$ are independent log-concave 
random variables, then their sum $V+W$ is also log-concave.
Equivalently, the convolution of any two log-concave sequences is log-concave.
\end{theorem}
A similar result was proved by Davenport and P\'{o}lya
\cite{davenport}, under the condition that the probability mass
functions of $V$ and $W$, rescaled by binomial coefficients, are
log-concave.

Notice that the sum of two independent and identically distributed
geometric random variables (both on the edge case) is a negative
binomial distribution (still log-concave, but no longer the edge case),
suggesting room for improvement.  In Section \ref{sec:proof}, we prove
the following main theorem, an extension of Theorem \ref{thm:hoggar}
to the case of dependent random variables.
\begin{theorem} \label{thm:main} If $V$ is a log-concave random variable and
random variables $V$, $W$ satisfy Condition \ref{cond:main} below then
the sum $V+W$ is also log-concave.
\end{theorem}
We give two distinct applications of Theorem \ref{thm:main} in the
remaining sections of the paper. 
Firstly, we consider the important property of 
log-concavity of certain combinatorial sequences.  It is often
desirable to show that such sequences
are unimodal (a sequence $c(\cdot)$ is unimodal if there exists an index
$j$ such that $c(i) \leq c(i+1)$ for $i \leq j-1$ and $c(i) \geq
c(i+1)$ for $i \geq j$).  Log-concavity is easily shown to imply
unimodality and, since unimodality is not preserved on convolution, 
log-concavity is
often easier to prove.  The motivation for Hoggar's
paper~\cite{hoggar} came from the coefficient sequences of chromatic
polynomials.  For a graph $G$, the chromatic polynomial $\pi_{G}(z) =
\sum_{i} c(i) z^i$ gives the number of ways to colour $G$ using exactly
$z$ colours.  Read~\cite{read} conjectured that $(|c(i)|, i \geq 0)$ is
a unimodal sequence and Welsh~\cite{welsh2} conjectured, more
strongly, that it is log-concave.  Hoggar proved
Theorem~\ref{thm:hoggar} in a partial attempt to resolve this
conjecture.  It remains unproven, but progress towards it is reviewed
in Brenti~\cite{brenti}.

In Section~\ref{sec:appl}, we show how a version of
Theorem~\ref{thm:main} (adapted to apply to sequences rather than just
probability mass functions) can be used to prove the log-concavity of
certain combinatorial sequences.  In particular, in Theorem~\ref{thm:stir2}
we use it to give a
simple proof of the log-concavity of the Stirling numbers of the
second kind (Condition~\ref{cond:main} is natural here and
straightforward to check). In Theorem \ref{thm:euler},
we prove log-concavity of the Eulerian numbers,
even though our Condition~\ref{cond:main} fails in this case, by a simple
adaptation of our methods.

The second application comes  in Section
\ref{sec:independent}, where for fixed log-concave $V$ we consider the set
of independent random variables $W$ such that $V+W$ is log-concave.
This means that we weaken one of Hoggar's assumptions, since we no
longer require $W$ to be log-concave.  We prove Lemma \ref{lem:equiv},
which gives simple conditions under which $X + \geo{p}$ is
log-concave, for independent $X$ and $\geo{p}$.
Further, we extend the arguments of 
Theorem \ref{thm:main} to prove Theorem
\ref{thm:order}, which states that if $X + \geo{p_1}$ is log-concave,
then so is $X + \geo{p_2}$ for $p_2 \geq p_1$. (This result is not
unexpected, but it cannot be stated as a simple corollary of Hoggar's 
Theorem). This result is relevant to the field of 
reliability theory, where we may
wish to establish the log-concavity (and hence bounds on the hazard rate) 
of the total lifetime of a collection of components, where some of the
components 
have a geometric lifetime with unknown (but bounded below) parameter $p$.
This restriction on the geometric parameter is very natural if we consider 
trying to infer $p$ from a collection of lifetimes while the trials are 
ongoing and some components ``haven't failed yet''.

Log-concavity is a property much used in econometrics, as described in
the survey by Bergstrom and Bagnoli~\cite{bergstrom}.  A previous
attempt to generalise Hoggar's theorem in this context, by Brais,
Martimort and Rochet~\cite{brais}, led to the statement of the
incorrect result that for independent $V$ and $W$
if either of $V$ or $W$ is log-concave then
their sum must be (see the discussion in Miravete~\cite{miravete}).
It is clear, however, that extensions of Hoggar's theorem may have
applications here.

A version of Theorem~\ref{thm:main} where $W$ must be a Bernoulli
random variable is proved in a different context by
Sagan~\cite{sagan1} (Theorem 1).  Our conditions reduce to his in this
special case.  Wang \cite{wang2} proves that linear transformations
based on binomial coefficients preserve log-concavity (see also
\cite{wang}). An anonymous referee brought it to our attention that,
at about the same time as the first draft of this paper, Wang and Yeh
had independently submitted a paper \cite{wang3} which contains some
similar technical results.  In particular, our Lemma \ref{lem:tech}
corresponds to their Lemma 2.1 and our Condition \ref{cond:main}
corresponds to their Lemma 2.3.
\section{Conditions and proof of Theorem \ref{thm:main}} \label{sec:proof}
We prove a generalized version of Theorem \ref{thm:main}, which
establishes both results about the log-concavity of the sum of
dependent random variables and results about the log-concavity of a
extended version of convolutions of sequences.

In general, given sequences $p_V(v)$ and $p_{V+W}(x)$, we will suppose
that there exists a two-dimensional array of coefficients
$p_{W|V}(w|v)$ such that
\begin{equation} \label{eq:newseq} 
p_{V+W}(x) = \sum p_V(v) p_{W|V}(x-v|v).\end{equation} In essence, the
$p_{W|V}$ act like conditional probabilities, but with no need for the
sequences to sum to 1.

We give here an example of a pair of dependent geometric random
variables whose sum is also geometric, motivating Condition
\ref{cond:main} below, under which we prove our extension of Hoggar's
theorem.
\begin{example} \label{ex:mot}
For some $p \in (0,1)$ and $\alpha \in (0,1)$, define the joint distribution
of $V$ and $W$ by 
$$
\pr( V = i, W = j) = \binom{i+j}{i} (1-p) p^{i+j} \alpha^i (1-\alpha)^j,
\mbox{ for $i,j \geq 0$}.$$
Using the identities $\sum_{i=0}^k \binom{k}{i} \alpha^i 
(1-\alpha)^{k-i} = 1$ and $\sum_{j=0}^{\infty} \binom{i+j}{i} t^j = 
(1-t)^{-i-1}$ (for $0 \leq t < 1$), we deduce that
$V+W$ is $\geo{p}$, $V$ is $\geo{\alpha p/(\alpha p + (1-p)}$ and
$W$ is $\geo{p/((1-\alpha) p + (1-p))}$.
The conditional probabilities are negative binomial with
\begin{equation} \label{eq:motcon}
p_{W|V}(j|i) = 
\pr( W = j | V = i) = \binom{i+j}{i} (1- p + \alpha p)^{i+1} 
(p(1-\alpha))^j.
\end{equation}
\end{example}
\begin{definition} Given coefficients $p_{W|V}$ and fixed $i$, 
define
$$ a^{(i)}_{r,s} = p_{W|V}(i-r|r)  p_{W|V}(i-s|s) -
p_{W|V}(i-r-1|r)  p_{W|V}(i-s+1|s).$$ \end{definition}
\begin{condition} \label{cond:main} For the
quantities $a^{(i)}_{r,s}$ defined above, 
we require that for all $0 \leq t \leq m \leq i$,
$$ \mbox{ (a) } \sum_{k = -t}^{t} a^{(i)}_{m+k,m-k} \geq 0 
\mbox{ and (b) } \sum_{k = -t-1}^{t} a^{(i)}_{m+k+1,m-k} \geq 0.
$$
\end{condition}
\begin{remark} \label{rem:s1s2}
Note that Condition \ref{cond:main} holds in Example \ref{ex:mot}. The
key is to observe that, in this case, for any given $i$,
$a^{(i)}_{r,s}$ is proportional to $\binom{i}{r} \binom{i}{s} -
\binom{i-1}{r} \binom{i+1}{s}$.  This means that for Condition 1(a),
as $t$ increases, the increment term $a^{(i)}_{m+t,m-t} +
a^{(i)}_{m-t,m+t}$ is proportional to $2\binom{i}{m+t} \binom{i}{m-t}
- \binom{i-1}{m+t} \binom{i+1}{m-t} - \binom{i+1}{m+t}
\binom{i-1}{m-t}$, which is positive for $t \leq T$ and negative for
$t > T$ (for some value $T$).  Hence, the partial sums,
$\sum_{k=-t}^{t} a_{m+k, m-k}^{(i)}$, form a sequence which increases
for $t \leq T$ and decreases thereafter.  Using the identity $\sum_{j}
\binom{a}{j} \binom{b}{r-j} = \binom{a+b}{r}$, we see that $ \sum_{k =
-m}^{m} a^{(i)}_{m+k,m-k} = 0$ and so the sequence of partial sums
must be non-negative for any $t$.
A similar argument holds for Condition 1(b).
\end{remark}
In order to prove Theorem \ref{thm:main} 
we require a technical lemma.
\begin{lemma}  \label{lem:tech}
Fix $l \geq m$ and suppose that $(c_j)$ is a sequence such that $C_n :=
\sum_{j=0}^n c_j \geq 0$ for all $0 \leq n \leq m$. For any
log-concave sequence $p$, and for any $0 \leq i \leq m$,
$$ \sum_{j=0}^i p(l+j) p(m-j) c_j \geq 0.$$ \end{lemma}
(This result is obvious if each $c_j \geq 0$, but the condition
$C_n \geq 0$ for all $0 \leq n \leq m$ is evidently weaker).

\begin{proof} 
We apply Abel's summation formula (``summation by parts'') to obtain
that 
\begin{align*} 
& \sum_{j=0}^i p(l+j) p(m-j) c_j = \sum_{j=0}^i p(l+j) p(m-j) (C_j - C_{j-1})  \\
& = \sum_{j=0}^i \big(p(l+j) p(m-j) - p(l+j+1) p(m-j-1)\big) C_j + C_i
p(l+i+1) p(m-i-1),
\end{align*}
where, by convention, $C_{-1} = 0$.  The log-concavity of $p$ implies
that $p(l+j) p(m-j) \geq p(l+j+1) p(m-j-1)$ for $j \geq 0$ and so,
since each $C_j \geq 0$, the result follows.
\end{proof}
\begin{proof}{\bf of Theorem \ref{thm:main} }
For any $i$, the sequence $p_{V+W}$ defined by (\ref{eq:motcon}) satisfies
\begin{eqnarray}
\lefteqn{p_{V+W}(i)^2 - p_{V+W}(i-1) p_{V+W}(i+1)} \nonumber \\
& = &
\sum_{j=0}^i \sum_{k=0}^{i+1} p_V(j) p_V(k)
\left\{ p_{W|V}(i-j|j)  p_{W|V}(i-k|k) \right. \nonumber \\
& & \hspace{4cm} -
\left. p_{W|V}(i-j-1|j)  p_{W|V}(i-k+1|k) \right\} \nonumber \\
& = & \sum_{j=0}^i \sum_{k=0}^{i+1} p_V(j) p_V(k) a^{(i)}_{j,k}. \label{eq:dec}
\end{eqnarray}
For simplicity, we decompose the above sum into three regions: 
(a) $\{ j \geq k \}$ (b) $\{ j = k-1 \}$
(c) $\{ j \leq k-2 \}$. 
We relabel the third region, using the new co-ordinates $(J,K) =
(k-1,j+1)$, which transforms it into $\{ J \geq K \}$. This enables us
to rewrite Equation (\ref{eq:dec}) as
\begin{align} 
\lefteqn{\sum_{j=0}^{i}  p_V(j) p_V(j+1) a^{(i)}_{j,j+1}
+ \sum_{i \geq j \geq k} \left( p_V(j) p_V(k) a^{(i)}_{j,k} +
p_V(j+1) p_V(k-1) a^{(i)}_{k-1,j+1} \right)} \label{eq:bigsum} \\
\lefteqn{= \sum_{m=0}^{i} \Bigg\{ \sum_{k=0}^{m} \left[ p_V(m+k) p_V(m-k)
a_{m+k, m-k}^{(i)} + p_V(m+k+1) p_V(m-k-1) a_{m-k-1, m+k+1}^{(i)}
\right] } \nonumber \\
& + p_V(m) p_V(m+1) a_{m, m+1}^{(i)} \nonumber \\
& +  \sum_{k=1}^{m+1} \left[ p_V(m+k) p_V(m-k+1) a_{m+k, m-k+1}^{(i)}
+ p_V(m+k+1) p_V(m-k) a_{m-k, m+k+1}^{(i)} \right] \Bigg\}.
\label{eq:bigsum2}
\end{align}
Here the second term in curly brackets corresponds to the first term
in Equation (\ref{eq:bigsum}). The first and third terms correspond
to the second term in Equation (\ref{eq:bigsum}), split according to
whether $r = j +k$ is even or odd. If $r$ is even then $m = r/2$; if
$r$ is odd then $m = (r-1)/2$. 
Consider now the first term in (\ref{eq:bigsum2}):
\begin{align}
\lefteqn{
\sum_{k=0}^{m} \left[ p_V(m+k) p_V(m-k) a_{m+k,m-k}^{(i)}
+ p_V(m+k+1) p_V(m-k-1) a_{m-k-1,m+k+1}^{(i)} \right]} \nonumber \\
& = \sum_{k=0}^{m} p_V(m+k) p_V(m-k) a_{m+k,m-k}^{(i)}
+ \sum_{k=0}^{m-1} p_V(m+k+1) p_V(m-k-1) a_{m-k-1,m+k+1}^{(i)} \nonumber \\
& = p_V(m)^2 a_{m,m}^{(i)} + 
\sum_{k=1}^{m} p_V(m+k) p_V(m-k) \left\{ a_{m+k,m-k}^{(i)}
+ a_{m-k,m+k}^{(i)} \right\}. \nonumber \\
& = \sum_{k=0}^{m} p_V(m+k) p_V(m-k) c_k, \label{eq:term1}
\end{align}
where $c_0 = a^{(i)}_{m,m}$ and
$c_k = a_{m+k,m-k}^{(i)} + a_{m-k,m+k}^{(i)}$ for $1 \leq k \leq m$.
Then Condition \ref{cond:main}(a) tells us that $\sum_{k=0}^t c_k \geq 0$
for all $0 \leq t \leq m$ and so, by Lemma \ref{lem:tech}
with $l = m$, $i=m$, Equation (\ref{eq:term1}) is positive. 
In the same way, we can 
show that sum of the second and third terms in (\ref{eq:bigsum2}) equals
\begin{equation}
\sum_{k=0}^{m} p_V(m+k+1) p_V(m-k) d_k, \label{eq:term2} 
\end{equation}
where $d_k = a_{m+k+1,m-k}^{(i)} + a_{m-k,m+k+1}^{(i)}$ for $0 \leq k \leq
m$.
Then Condition \ref{cond:main}(b) tells us that $\sum_{k=0}^t d_k \geq 0$
for all $0 \leq t \leq m$ and so, by Lemma \ref{lem:tech}
with $l = m+1$, $i=m$, Equation (\ref{eq:term2}) is positive.

Hence, $p_{V+W}(i)^2 - p_{V+W}(i-1) p_{V+W}(i+1) \geq 0$ for all $i \geq 0$.
\end{proof}

We now identify which properties of independent random variables allow
Hoggar's result, Theorem \ref{thm:hoggar}, to be proved. In particular,
Condition \ref{cond:main} holds under Hoggar's assumptions, so Theorem
\ref{thm:main} implies Theorem \ref{thm:hoggar}.
\begin{remark} \label{rem:indprop}
For independent random variables
$V$ and $W$, the $a^{(i)}_{j,k}$ have the following
properties: (a) for all $j$, $a^{(i)}_{j,j+1} \equiv 0$ (b) for $j \geq k$, 
$a^{(i)}_{k-1,j+1} = 
- a^{(i)}_{j,k}$ (c) if $W$ is log-concave then $a^{(i)}_{j,k} \geq 0$ for $j \geq k$.

Hence, if we fix $i$ and define $c_j$ and $d_j$ as in the proof of Theorem
\ref{thm:main}, for $V$ and $W$ independent and log-concave we can write:
\begin{eqnarray*}
\sum_{j=0}^t c_j & = & a^{(i)}_{m,m} + \sum_{j=1}^t a^{(i)}_{m+j,m-j} -
\sum_{j=1}^t a^{(i)}_{m+j-1,m-j+1} \\
& = & a^{(i)}_{m,m} + \sum_{j=1}^t a^{(i)}_{m+j,m-j} -
\sum_{j=0}^{t-1} a^{(i)}_{m+j,m-j} =  a^{(i)}_{m+t,m-t} \geq 0.
\end{eqnarray*}
Similarly:
\begin{eqnarray*}
\sum_{j=0}^t d_j & = & \sum_{j=0}^t a^{(i)}_{m+1+j,m-j} -
\sum_{j=0}^t a^{(i)}_{m+j,m-j+1} \\
& = & a^{(i)}_{m+1+t,m-t} - a^{(i)}_{m,m+1} = a^{(i)}_{m+1+t,m-t} \geq 0.
\end{eqnarray*}
Thus, Condition \ref{cond:main} holds for independent and log-concave
$V$ and $W$.
\end{remark}
\section{Log-concavity of combinatorial sequences} \label{sec:appl}
Many common combinatorial sequences, such as the binomial coefficients
$(\binom{n}{k}, 0 \leq k \leq n)$, are log-concave in $k$.  
A key result for proving the log-concavity of combinatorial sequences
is the following.
\begin{lemma} \label{lem:realroots}
If all the roots of the generating function
$p(x) = c(0) x^n + c(1) x^{n-1} + c(2) x^{n-2} + \cdots + c(n)$
are real and negative, then its coefficients $c(k)$ form a log-concave 
sequence.
\end{lemma}

The proof can be regarded as an application of
Theorem~\ref{thm:hoggar}.  See Niculescu~\cite{niculescu} for a
review of this and similar properties.

Section
4.5 of Wilf~\cite{wilf} contains proofs of
log-concavity for the binomial coefficients and the Stirling numbers
of the first and second kinds. He shows that the generating function
of each sequence has real and negative roots, and then appeals to
Lemma \ref{lem:realroots} to deduce log-concavity of the sequence.

For the first two examples he considers, this is relatively
straightforward.  That is, in Corollary 4.5.1 he uses the fact that
the binomial coefficients $\binom{n}{k}$ have generating function
$(1+x)^n$, and in Corollary 4.5.2 he uses the fact that the Stirling
numbers of the first kind have generating function $(x+1)(x+2) \ldots
(x + n-1)$.

However, the fact that these generating functions factorise into
products of independent terms indicate that both these examples can
equally well be dealt with via Hoggar's Theorem. As remarked in
Example \ref{ex:first}.\ref{ex:bernsum}, the binomial distribution is
log-concave.  Further, the fact that the Stirling numbers of the first
kind satisfy the recurrence relation $S_1(n+1,k) = S_1(n,k-1) + n S_1(n,k)$
means that we can write $ S_1(n,k) = n! \pr( X_1 + \ldots + X_n = k)$
(the normalization is irrelevant here), where $X_i$ are independent
with $\pr(X_i = 0) = (i-1)/i$ and $\pr(X_i = 1) = 1/i$, so the result
derived by Wilf is simply another case of Example
\ref{ex:first}.\ref{ex:bernsum}.

Indeed, Hoggar's Theorem can prove log-concavity in cases where the
method described by Wilf fails.  For fixed $n$, let $b(n,k)$ be the
number of permutations of $n$ objects with exactly $k$ inversions
(that is, the permutation $\sigma$ has exactly $k$ pairs $i < j$ such
that $\sigma(i) > \sigma(j)$).  Equation (4.5.7) of \cite{wilf} shows
that the generating function of $b(n,k)$ takes the form
$$ p(x) = \prod_{i=1}^{n-1} (1 + x + \ldots + x^{i}).$$ This means that the
sequence $b(n,k)$ can be expressed as the convolution of independent
log-concave sequences of the form $(1, 1, \ldots, 1)$, and so is
log-concave. However, for $n \geq 3$, 
the generating function $p(x)$ has imaginary roots, so
Lemma~\ref{lem:realroots} does not apply.

The Stirling numbers of the second kind provide an example where Wilf's 
proof is relatively complicated and where Hoggar's result does not apply.
However, our Theorem \ref{thm:main} gives a simple proof
in this case, using only the property that
$S_2(n,k) = S_2(n-1,k-1) + k S_2(n-1,k)$. In fact, this property is used
by Sagan \cite{sagan1} in his proof of the same result; however, his techniques
only apply in the case where $W$ is a $\{0,1\}$-valued random variable,
whereas our Theorem \ref{thm:main} will prove log-concavity of a more
general family of combinatorial sequences.
\begin{theorem} \label{thm:stir2}
For given $n$, the Stirling numbers
of the second kind, 
$S_2(n,k)$ form a log-concave sequence in $k$. \end{theorem}
\begin{proof} 
In the notation of Equation (\ref{eq:newseq}), if $p_V(i) = S_2(n-1,i)$
then $p_{V+W}(i) = S_2(n,i)$, where for each $j \leq n-1$,
$$ p_{W|V}(x|j) = \left\{ 
\begin{array}{ll} 
j & \mbox{ if $x = 0$,} \\
1 & \mbox{ if $x = 1$,} \\
0 & \mbox{ otherwise,}
\end{array} \right. $$
(so that $V$ and $W$ are not ``independent'').  Given $i$, for $j \geq
k$, the only non-zero $a_{j,k}^{(i)}$ are 
\begin{center}
\begin{tabular}{c||c|c|c}
$a^{(i)}_{j,k}$ & $k=i-1$ & $k=i$ & $k=i+1$ \\
\hline \hline 
$j=i-2$ & $0$ & $-1$ & $-(i+1)$ \\
$j=i-1$ & $1$ & $i - (i-1)$ & $-(i+1)(i-1)$ \\
$j=i$   & $i$ & $i^2$ & $0$ \\
\end{tabular}
\end{center}
Hence, the only non-zero sequences of partial sums in Condition 1(a)
are the case $m = i-1$, where $\left(\sum_{k=-t}^{t}
a^{(i)}_{m+k,m-k}, t \geq 0\right)$ is $(1, 0, 0, \ldots)$, and the
case $m=i$, where $\left(\sum_{k=-t}^{t} a^{(i)}_{m+k,m-k}, t \geq
0\right)$ is $(i^2, 1, 1, \ldots)$.  In Condition 1(b), the only
non-zero sequence of partial sums occurs for $m = i-1$, where
$\left(\sum_{k=-t-1}^{t} a^{(i)}_{m+k+1,m-k}, t \geq 0\right)$ is $(i,
i+1, 0, 0, \ldots)$.  These are all positive, so Condition
\ref{cond:main} holds, as required.
\end{proof}
This result was first proved by Lieb \cite{lieblc};
see Wang and Yeh~\cite{wang} for an alternative approach.

The proof of Theorem \ref{thm:stir2} 
can also be used to show that the $q$-Stirling numbers of
the second kind, $(S_2^{(q)}(n,k), 0 \leq k \leq n)$, form a log-concave
sequence for $0 \leq q \leq 1$ (see Sagan~\cite{sagan2} for the
definition of the $q$-Stirling numbers and a proof of this fact for
all $q \geq 0$).

Indeed, our methods can sometimes be adapted to prove log-concavity of 
combinatorial sequences even when Condition \ref{cond:main} and
Theorem~1 of \cite{sagan1} fail. 
For 
example, consider the Eulerian numbers $E(n,k)$, the number of permutations of
$n$ objects with exactly $k$ ascents 
(that is, the permutation $\sigma$ has exactly $k$ places $j$ such
that $\sigma(j) < \sigma(j+1)$). 
We prove the following well-known result 
(see, for example, 
Gasharov \cite{gasharov} or B\'{o}na and Ehrenborg \cite{bona} for
alternative proofs).
\begin{theorem} \label{thm:euler}
For given $n$, the Eulerian numbers
$E(n,k)$ form a log-concave sequence in $k$. \end{theorem}
\begin{proof} 
The recurrence 
relation $E(n,k) = (k+1) E(n-1,k) + (n-k) E(n-1,k-1)$, implies that
$$ p_{W|V}(x|j) = \left\{ 
\begin{array}{ll} 
j+1 & \mbox{ if $x = 0$,} \\
n-j & \mbox{ if $x = 1$,} \\
0 & \mbox{ otherwise.}
\end{array} \right. $$
Given $i$, for $j \geq
k$, the only non-zero $a_{j,k}^{(i)}$ are 
\begin{center}
\begin{tabular}{c||c|c|c}
$a^{(i)}_{j,k}$ & $k=i-1$ & $k=i$ & $k=i+1$ \\
\hline \hline 
$j=i-2$ & $0$ & $-(n-i)(n-i+2)$ & $-(i+2)(n-i+2)$ \\
$j=i-1$ & $(n-i+1)^2 $ & $n+1$ & $-i(i+2)$ \\
$j=i$   & $(i+1)(n-i+1)$ & $(i+1)^2$ & $0$ \\
\end{tabular}
\end{center}
For $m = i-1$, Condition 1(b) fails.
However, since we are attempting to show that 
$\sum_{j,k} a_{j,k}^{(i)} p_V(j) p_V(k)$ is positive, it will be sufficient
to show that this sum minus a positive quantity is positive. If we define
\begin{eqnarray*}
\sum_{j=0}^i \sum_{k=0}^{i+1} \widetilde{a}_{j,k}^{(i)} p_V(j) p_V(k)
=  \sum_{j=0}^i \sum_{k=0}^{i+1}  a_{j,k}^{(i)} p_V(j) p_V(k)
- (p_V(i-1) - p_V(i))^2,
\end{eqnarray*}
then the new coefficients are
\begin{center}
\begin{tabular}{c||c|c|c}
$\widetilde{a}^{(i)}_{j,k}$ & $k=i-1$ & $k=i$ & $k=i+1$ \\
\hline \hline 
$j=i-2$ & $0$ & $-(n-i)(n-i+2)$ & $-(i+2)(n-i+2)$ \\
$j=i-1$ & $(n-i+1)^2-1 $ & $(n+1)+1$ & $-i(i+2)$ \\
$j=i$   & $(i+1)(n-i+1)+1$ & $(i+1)^2-1$ & $0$ \\
\end{tabular}
\end{center}
These new coefficients $\widetilde{a}^{(i)}_{j,k}$
do satisfy Condition \ref{cond:main}, so 
we deduce that the 
original difference $p_{V+W}(i)^2 - p_{V+W}(i-1) p_{V+W}(i+1)$
is positive.
\end{proof}
\section{Development for the independent case} \label{sec:independent}
We now use a version of Theorem \ref{thm:main} to give results about 
when the sum of independent random variables $V$ and $W$
is log-concave, even though $V$ and $W$ need not be. 

To give a concrete example of an independent pair of random variables
$V$ and $W$ such that $V+W \in LC$ but $W \not \in LC$, let $W$ have
probability mass function $ p_W(0) = 5/8$, $p_W(1) = 1/4$, $p_W(2) =
1/8$.  Then $p_W(1)^2 - p_W(0) p_W(2) = -1/64$ and so $W$ is not
log-concave.  However, $p_W(1)/p_W(0) = 2/5$ and $p_W(2)/p_W(1) = 1/2$
and so, by Lemma \ref{lem:equiv} below, $W + \geo{p}$ is log-concave
for any $p \geq 1/2$.
\begin{definition} For each random variable $V$, define the set
$$ \C_V = \{ \text{ $W$ independent of $V$}: \text{$V+W$  is log-concave } \}.$$
Define a partial order on random variables by $V_1 \preceq V_2$ if and only
if $\C_{V_1} \subseteq \C_{V_2}.$ 
Write $LC$ for the set of log-concave variables. \end{definition}
We first give some simple properties of the sets $\C_V$, before going on
to consider properties of the sets $\C_{\geo{p}}$.

\begin{proposition} For the sets $\C_V$ defined above:
\begin{enumerate}
\item{If $W \in \C_V$ then $V \in \C_W$.} 
\item{For each $V$, $\C_V$ is a closed set.}
\item{The $\C_V$ are partial order ideals with respect to $\preceq$.
That is, given $W_1 \preceq W_2$, if $W_1 \in \C_V$ then $W_2 \in \C_V$.}
\item{$ \begin{displaystyle} \bigcap_{V \in LC} \C_V = LC \end{displaystyle}$.
} 
\item{If $V_2 = V_1 + U$, where $V_1$ and $U$ are independent and $U$ is log-concave,
then $V_1 \preceq V_2$. }
\end{enumerate}
\end{proposition}
\begin{proof}
The first result is trivial. For each $i$, the set of $W$ with
$p_{V+W}(i)^2 \geq p_{V+W}(i-1) p_{V+W}(i+1)$ is closed, since it can
be expressed as the inverse image of the closed set $[0,\infty)$ under
a continuous map which depends only on $p_{W}(0), \ldots, p_W(i+1)$.
This means that $\C_V$ is a countable intersection of closed sets and
so Part 2 follows.

Part 3 is an application of Part 1, since 
\[
 W_1 \in \C_V \Rightarrow V \in \C_{W_1} \Rightarrow V \in \C_{W_2}
\Rightarrow W_2 \in C_V.
\]
To prove Part 4, note that Theorem \ref{thm:hoggar} entails that for
all $V \in LC$, $LC \subseteq \C_V$, so $LC \subseteq \bigcap_{V \in
LC} \C_V$.  Further, since for $V \equiv 0$ we have $\C_V= LC$, this
set inclusion must be an equality.

For any $W \in \C_{V_1}$, the random variable
$V_1 + W$ is log-concave and so by Theorem \ref{thm:hoggar} 
$(V_1 + W) + U = (V_2 + W)$ is log-concave. Hence
$W \in \C_{V_2}$ and Part 5 follows. \end{proof}
\begin{lemma} \label{lem:equiv}
For a random variable $X$ with probability mass function $p_X(i)$:
\begin{enumerate}
\item{If $p_X(i+1)/p_X(i) \leq p$ for all $i \geq 1$, then $X \in 
\C_{\geo{p}}$.}
\item{If there is a gap in the support of $X$ (that is
for some $i$, $p_X(i) = 0$, but $p_X(i-1) > 0$ and $p_X(i+1) > 0$) then
$X \notin \C_{\geo{p}}$ for any $p$.}
\end{enumerate}
\end{lemma}
\begin{proof}
Write $q_p(i)$ for the probability mass function of $X + \geo{p}$.
Note that since 
\begin{equation} \label{eq:recrel}
q_p(i) = \sum_{j=0}^i (1-p)p^{i-j} p_X(j) = (1-p)p_X(i)
+ p q_p(i-1), \end{equation} we have
\begin{eqnarray}
\lefteqn{q_p(i)^2 - q_p(i+1) q_p(i-1)} \nonumber \\ 
& = & q_p(i) \{ (1-p)p_X(i) + p q_p(i-1) \}
- q_p(i-1) \{ (1-p)p_X(i+1) + p q_p(i) \} \nonumber \\
& = & (1-p) \{ q_p(i) p_X(i) - q_p(i-1) p_X(i+1) \}. \label{eq:decomp}
\end{eqnarray}
This means that the sum $X+ \geo{p}$ is log-concave if and only if for
all $i \geq 1$,
\begin{equation} \label{eq:equiv} q_p(i) p_X(i) - q_p(i-1)
p_X(i+1) \geq 0. \end{equation}

By (\ref{eq:recrel}), we have that $q_p(i) \geq p q_p(i-1)$ and so if
$p_X(i+1)/p_X(i) \leq p$ for all $i \geq 1$, then
$$ \frac{p_X(i+1)}{p_X(i)} \leq p \leq
\frac{q_p(i)}{q_p(i-1)},$$ and so $X+ \geo{p}$ is log-concave.

If $p_X(i-1) > 0$ then $q_p(i-1) > 0$. Hence, if $p_X(i) = 0$ and
$p_X(i+1) > 0$ then the inequality (\ref{eq:equiv}) fails to hold for
this particular $i$, so $X + \geo{p}$ is not log-concave.
\end{proof}
We now use arguments suggested by the proof of Theorem \ref{thm:main} to 
show that the 
property that $X +\geo{p}$ is log-concave holds monotonically in $p$. The
key, as in Equation (\ref{eq:dec}), is to expand 
$ q_{p_2}(i)^2 - q_{p_2}(i-1) q_{p_2}(i+1)$ as a sum of terms
of the form $q_{p_1}(j) q_{p_1}(k) a_{j,k}^{(i)}$. This is achieved
in Equation (\ref{eq:tosort}). 
\begin{theorem} \label{thm:order}
The partial order $\preceq$ induces a total order on
the set of geometric random variables, that is:
$$ \geo{p_1} \preceq \geo{p_2} \mbox{ if and only if } p_1 \leq p_2.$$
\end{theorem}
\begin{proof}
With the notation of Lemma \ref{lem:equiv}, if 
$p_1 \leq p_2$, there exist positive constants $b_r$, where 
$b_0 = (1-p_2)/(1-p_1)$,
$b_r = \left((p_2 - p_1)(1-p_2)/(1-p_1) \right) p_2^{r-1}$, such that
\begin{equation} \label{eq:qp2} q_{p_2}(i) = 
 \sum_{r=0}^i b_r q_{ p_1}(i-r).\end{equation}
This follows using a standard coupling argument, or since
\begin{eqnarray*}
\frac{ (p_2 - p_1)(1-p_2)}{(1-p_1)}
\sum_{r=0}^i p_2^{r-1} q_{p_1}(i-r) 
& = & \frac{ (p_2 - p_1)(1-p_2)}{p_2}
\sum_{j=0}^i q(j) \sum_{j=0}^{i-j} p_2^r p_1^{i-r-j} \\ 
& = & \frac{ (1-p_2)}{p_2}
\sum_{j=0}^i q(j) \left( p_2^{i-j+1} - p_1^{i-j+1} \right) \\ 
& = & q_{p_2}(i) - \frac{ p_1(1-p_2)}{p_2(1-p_1)} q_{p_1}(i).
\end{eqnarray*}
Now, we can combine Equations (\ref{eq:recrel}),
(\ref{eq:decomp}) and (\ref{eq:qp2}) to obtain:
\begin{eqnarray} 
\lefteqn{q_{p_2}(i)^2 - q_{p_2}(i+1) q_{p_2}(i-1)} \nonumber \\ 
& = & (1-p_2) \{ q_{p_2}(i) p_X(i) - q_{p_2}(i-1) p_X(i+1) \} \nonumber \\
& = & \frac{(1-p_2)}{(1-p_1)} \sum_{r=0}^i b_r  \biggl\{ q_{p_1}(i-r) 
\big(q_{p_1}(i) - p_1 q_{p_1}(i-1) \big) \nonumber \\
& & \hspace*{4cm} 
- q_{p_1}(i-r-1) \big(q_{p_1}(i+1) - p_1 q_{p_1}(i) \big) \biggr\}.  
\label{eq:tosort}
\end{eqnarray}
Now, each term in curly brackets is positive, since
$q_{p_1}(i) - p_1 q_{p_1}(i-1) \geq 0$, so 
by log-concavity of $q_{p_1}$, for all $r \geq 0$:
$$
\frac{ q_{p_1}(i-r)}{ q_{p_1}(i-r-1)} \geq \ldots \geq 
\frac{ q_{p_1}(i)}{ q_{p_1}(i-1)} \geq 
\frac{ q_{p_1}(i+1) - p_1 q_{p_1}(i)}{ q_{p_1}(i) - p_1 
q_{p_1}(i-1)},$$
and so $q_{p_2}(i)^2 - q_{p_2}(i+1) q_{p_2}(i-1) \geq 0$.
\end{proof}
\section*{Acknowledgments}
We are grateful to anonymous referees whose comments led to considerable
improvements in the presentation of this paper.
\bibliography{../../bibliography/phd}

\end{document}